\documentclass[oneside,american,english]{amsart}
\usepackage[T1]{fontenc}
\usepackage[latin9]{inputenc}
\usepackage{amsthm}
\usepackage{amstext}
\usepackage{amssymb}
\usepackage{esint}

\makeatletter
\numberwithin{equation}{section}
\numberwithin{figure}{section}
\theoremstyle{plain}
\newtheorem{thm}{\protect\theoremname}
  \theoremstyle{definition}
  \newtheorem{defn}[thm]{\protect\definitionname}
  \theoremstyle{plain}
  \newtheorem*{thm*}{\protect\theoremname}
  \theoremstyle{plain}
  \newtheorem{prop}[thm]{\protect\propositionname}
  \theoremstyle{plain}
  \newtheorem{lem}[thm]{\protect\lemmaname}

\usepackage{ae,aecompl} 

\usepackage{enumitem}

\renewenvironment{enumerate}{\begin{oldenumerate}[topsep=0pt]}{\end{oldenumerate}}

\makeatother

\usepackage{babel}
  \addto\captionsamerican{\renewcommand{\definitionname}{Definition}}
  \addto\captionsamerican{\renewcommand{\lemmaname}{Lemma}}
  \addto\captionsamerican{\renewcommand{\propositionname}{Proposition}}
  \addto\captionsamerican{\renewcommand{\theoremname}{Theorem}}
  \addto\captionsenglish{\renewcommand{\definitionname}{Definition}}
  \addto\captionsenglish{\renewcommand{\lemmaname}{Lemma}}
  \addto\captionsenglish{\renewcommand{\propositionname}{Proposition}}
  \addto\captionsenglish{\renewcommand{\theoremname}{Theorem}}
  \providecommand{\definitionname}{Definition}
  \providecommand{\lemmaname}{Lemma}
  \providecommand{\propositionname}{Proposition}
  \providecommand{\theoremname}{Theorem}
\providecommand{\theoremname}{Theorem}

\begin{document}
\selectlanguage{american}%
\global\long\def\epsilon{\varepsilon}

\global\long\def\phi{\varphi}

\global\long\def\R{\mathbb{R}}

\global\long\def\N{\mathbb{N}}

\global\long\def\o{\mathbf{1}}

\global\long\def\L{\mathcal{L}}

\global\long\def\kc{\mathcal{K}_{c}^{n}}

\global\long\def\vol{\text{Vol}}

\global\long\def\lc{\operatorname{LC}\left(\R^{n}\right)}

\global\long\def\lco{\operatorname{LC}_{0}\left(\R^{n}\right)}

\global\long\def\qc{\operatorname{QC}\left(\R^{n}\right)}

\global\long\def\qco{\operatorname{QC}_{0}\left(\R^{n}\right)}

\global\long\def\cvx{\operatorname{Cvx}\left(\R^{n}\right)}

\global\long\def\KL{\underline{K}}

\global\long\def\KU{\overline{K}}

\global\long\def\Ca{\operatorname{C}_{\alpha}\left(\R^{n}\right)}

\global\long\def\C#1{\operatorname{C}_{#1}\left(\R^{n}\right)}

\global\long\def\basea{\operatorname{base}_{\alpha}}

\global\long\def\base#1{\operatorname{base}_{#1}}

\global\long\def\GL{\operatorname{GL}}

\selectlanguage{english}%

\title{$\alpha$-concave functions and a functional extension of mixed volumes}

\author{Vitali Milman}

\address{School of Mathematical Sciences \\
Tel Aviv University \\
Tel Aviv 69978, Israel}

\email{milman@post.tau.ac.il}

\thanks{Partially supported by the Minkowski Center at the University of
Tel Aviv, by ISF grant 387/09 and by BSF grant 2006079.}

\author{Liran Rotem}

\address{School of Mathematical Sciences \\
Tel Aviv University \\
Tel Aviv 69978, Israel}

\email{liranro1@post.tau.ac.il}
\begin{abstract}
In this paper we define an addition operation on the class of quasi-concave
functions. While the new operation is similar to the well-known sup-convolution,
it has the property that it polarizes the Lebesgue integral. This
allows us to define mixed integrals, which are the functional analogs
of the classic mixed volumes.

We extend various classic inequalities, such as the Brunn-Minkowski
and the Alexandrov-Fenchel inequality, to the functional setting.
For general quasi-concave functions, this is done by restating those
results in the language of rearrangement inequalities. Restricting
ourselves to log-concave functions, we prove generalizations of the
Alexandrov inequalities in a more familiar form. 
\end{abstract}


\begin{abstract}
Mixed volumes, which are the polarization of volume with respect to
the Minkowski addition, are fundamental objects in convexity. In this
note we announce the construction of mixed integrals, which are functional
analogs of mixed volumes. We build a natural addition operation $\oplus$
on the class of quasi-concave functions, such that every class of
$\alpha$-concave functions is closed under $\oplus$. We then define
the mixed integrals, which are the polarization of the integral with
respect to $\oplus$. 

We proceed to discuss the extension of various classic inequalities
to the functional setting. For general quasi-concave functions, this
is done by restating those results in the language of rearrangement
inequalities. Restricting ourselves to $\alpha$-concave functions,
we state a generalization of the Alexandrov inequalities in their
more familiar form. 
\end{abstract}
\maketitle

\section{$\alpha$-concave functions}

Let us begin by introducing our main objects of study:
\begin{defn}
\label{def:a-concave}Fix $-\infty\le\alpha\le\infty$. We say that
a function $f:\R^{n}\to[0,\infty)$ is $\alpha$-concave if $f$ is
supported on some convex set $\Omega$, and for every $x,y\in\Omega$
and $0\le\lambda\le1$ we have
\[
f\left(\lambda x+(1-\lambda)y\right)\ge\left[\lambda f(x)^{\alpha}+\left(1-\lambda\right)f(y)^{\alpha}\right]^{\frac{1}{\alpha}}.
\]
For simplicity, we will always assume that $f$ is upper semicontinuous,
$\max_{x\in\R^{n}}f(x)=1,$ and $f(x)\to0$ as $\left|x\right|\to\infty$.
The class of all such $\alpha$-concave functions will be denoted
by $\Ca$. 
\end{defn}
In the above definition, we follow the convention set by Brascamp
and Lieb (\cite{Brascamp1976}), but the notion of $\alpha$-concavity
may be traced back to Avriel (\cite{Avriel1972}) and Borell (\cite{Borell1974},
\cite{Borell1975}). Discussions of $\alpha$-concave functions from
a geometric point of view may be found, e.g., in \cite{Bobkov2009}
and \cite{Rotem2012a}.

In the cases $\alpha=-\infty,0,\infty$ we understand Definition \ref{def:a-concave}
in the limit sense. For example, $f\in\C{\infty}$ if $f$ is supported
on some convex set $\Omega$, and 
\[
f\left(\lambda x+(1-\lambda)y\right)\ge\max\left\{ f(x),f(y)\right\} 
\]
 for all $x,y\in\Omega$ and $0\le\lambda\le1$. Of course, this just
means that $f$ is constant on $\Omega$. In other words, we have
a natural correspondence between $\C{\infty}$ to the class $\kc$
of compact, convex sets in $\R^{n}$: every function $f\in\C{\infty}$
is of the form 
\[
f(x)=\o_{K}(x)=\begin{cases}
1 & x\in K\\
0 & \text{otherwise,}
\end{cases}
\]
 for some $K\in\kc$. 

Notice that if $\alpha_{1}<\alpha_{2}$ then $\C{\alpha_{1}}\supset C_{\alpha_{2}}\left(\R^{n}\right)$
(see \cite{Brascamp1976}) . Therefore, we can view the class $\Ca$
for $\alpha<\infty$ as an extension of the class $\kc$ of convex
sets. Our main goal in this note is to extend the geometric notion
of mixed volumes from $\kc$ to the different classes of $\alpha$-concave
functions.

For convenience, we will restrict ourselves to the case $-\infty\le\alpha\le0$.
For $-\infty<\alpha<0$, it is easy to see that $f$ is $\alpha$-concave
if and only if $f^{\alpha}$ is a convex function on $\R^{n}$. The
cases $\alpha=0,-\infty$ are important, and deserve a special name:
\begin{defn}

\begin{enumerate}
\item A $0$-concave function is called log-concave. These are the functions
$f:\R^{n}\to[0,\infty)$ such that 
\[
f\left(\lambda x+(1-\lambda)y\right)\ge f(x)^{\lambda}f(y)^{1-\lambda}
\]
 for all $x,y\in\R^{n}$ and $0\le\lambda\le1$. We will usually write
$\lc$ instead of $\C 0$.
\item A $(-\infty)$-concave function is called quasi-concave. These are
the functions $f:\R^{n}\to[0,\infty)$ such that 
\[
f\left(\lambda x+(1-\lambda)y\right)\ge\min\left\{ f(x),f(y)\right\} 
\]
 for all $x,y\in\R^{n}$ and $0\le\lambda\le1$. We will usually write
$\qc$ instead of $\C{-\infty}$.
\end{enumerate}
\end{defn}
We will now see that if $-\infty<\alpha\le0$, there is a natural
correspondence between $\Ca$ and convex functions on $\R^{n}$. Since
we only care about negative values of $\alpha$, it will sometimes
be convenient to use the parameter $\beta=-\frac{1}{\alpha}$. The
following definition appeared in \cite{Rotem2012a}:
\begin{defn}
The convex base of a function $f\in\Ca$ is 
\[
\basea\left(f\right)=\frac{1-f^{\alpha}}{\alpha}.
\]
 Put differently, $\phi=\basea\left(f\right)$ is the unique convex
function such that 
\[
f=\left(1+\frac{\phi}{\beta}\right)^{-\beta}.
\]
 In the limiting case $\alpha=0$ we define $\base 0\left(f\right)=-\log f$. 
\end{defn}
By our assumptions on $f$, the function $\phi=\basea\left(f\right)$
is convex, lower semicontinuous, with $\min\phi=0$ and $\phi(x)\to\infty$
as $\left|x\right|\to\infty$. We will denote this class of convex
functions by $\cvx$ (this is not an entirely standard notation),
and notice that the map $\basea$ is a bijection between $\Ca$ and
$\cvx$. It follows immediately, for example, that every function
$f\in\Ca$ is continuous on its support, because the same is true
for convex functions.

In the case $\alpha=-\infty$, we have no such correspondence. It
is therefore not surprising that it possible to construct quasi-concave
functions which are not continuous on their support. Indeed, fix convex
sets $K_{1}\subset K_{2}$ and define $f=\o_{K_{1}}+\o_{K_{2}}$. 

Remember that if $f\in\Ca$, then $f\in\C{\alpha^{\prime}}$ for every
$\alpha'<\alpha$. However, in general we have $\basea\left(f\right)\ne\base{\alpha^{\prime}}(f)$,
so the base depends on the class we choose to work in, and not only
on our function $f$. However, in the specific case $f=\o_{K}$ for
some convex set $K$, we have 
\[
\basea(f)=\o_{K}^{\infty}=\begin{cases}
0 & x\in K\\
\infty & \text{otherwise}
\end{cases}
\]
 for every value of $\alpha$. 

On $\cvx$ there is a natural addition operation, known as inf-convolution:
\begin{defn}
\label{def:inf-conv}For $\phi,\psi\in\cvx$ we define their inf-convolution
to be
\[
\left(\phi\square\psi\right)(x)=\inf_{y+z=x}\left[\phi(y)+\psi(z)\right].
\]
 Similarly, if $\phi\in\cvx$ and $\lambda>0$ we will define
\[
\left(\lambda\cdot\phi\right)(x)=\lambda\phi\left(\frac{x}{\lambda}\right).
\]

\end{defn}
The definition of $\lambda\cdot\phi$ was chosen to have $2\cdot\phi=\phi\square\phi$,
as one easily verifies. It is also easy to see that we have commutativity,
associativity and distributivity.

We will not explain the exact sense in which these operations are
natural, and instead refer the reader to the first section of \cite{Rotem2012a}.
We will note, however, that Definition \ref{def:inf-conv} extends
the classical operations on convex bodies: If $K_{1},K_{2}\in\kc$
and $\lambda>0$ then 
\[
\left(\lambda\cdot\o_{K_{1}}^{\infty}\right)\square\o_{K_{2}}^{\infty}=\o_{\lambda K_{1}+K_{2}}^{\infty}.
\]
 Here $+$ is the Minkowski sum of convex bodies, defined by 
\[
K_{1}+K_{2}=\left\{ x+y:\ x\in K_{1},\ y\in K_{2}\right\} ,
\]
 and $\lambda K$ is defined by 
\[
\lambda K=\left\{ \lambda x:\ x\in K\right\} .
\]

We will now define addition on $\Ca$, using the established correspondence
between $\Ca$ and $\cvx$:
\begin{defn}
\label{def:a-sum}Fix $-\infty<\alpha\le0$. Then: 
\begin{enumerate}
\item For $f,g\in\Ca$ we define their $\alpha$-sum $f\star_{\alpha}g$
by the relation 
\[
\basea\left(f\star_{\alpha}g\right)=\left(\basea f\right)\square\left(\basea g\right).
\]

\item For $f\in\Ca$ and $\lambda>0$ we define $\lambda\cdot_{\alpha}f$
via the relation 
\[
\base{\alpha}\left(\lambda\cdot_{\alpha}f\right)=\lambda\cdot\base{\alpha}f.
\]

\end{enumerate}
\end{defn}
Again, the definition of $\alpha$-sum depends on $\alpha$, and not
only on $f$ and $g$: If $\alpha'<\alpha$ and $f,g\in\Ca$, then
in general we have $f\star_{\alpha}g\ne f\star_{\alpha^{\prime}}g$.
However, for indicators of convex sets we have 
\[
\o_{K_{1}}\star_{\alpha}\o_{K_{2}}=\o_{K_{1}+K_{2}}
\]
 for all $\alpha$. 

The definition of $\alpha$-sum may be written down explicitly, without
referring to the convex bases. For $-\infty<\alpha<0$ we have 
\begin{equation}
\left(f\star_{\alpha}g\right)(x)=\sup_{y+z=x}\left(f(y)^{\alpha}+g(z)^{\alpha}-1\right)^{\frac{1}{\alpha}},\label{eq:a-sum}
\end{equation}
 and for $\alpha=0$ we get the limiting case 
\[
\left(f\star_{0}g\right)(x)=\sup_{y+z=x}f(y)g(z).
\]
 The operation $\star_{0}$ on $\lc$ is known as Asplund-sum, or
sup-convolution (see, e.g., \cite{Klartag2005}). 

For $\alpha=-\infty$ we cannot define $f\star g$ using the same
approach as Definition \ref{def:a-sum}, because we do not have the
notion of a base for quasi-concave functions. However, we may use
equation \ref{eq:a-sum}, and the fact that for every $0<u,v\le1$
we have 
\[
\lim_{\alpha\to-\infty}\left[u^{\alpha}+v^{\alpha}-1\right]^{\frac{1}{\alpha}}=\min\left\{ u,v\right\} .
\]
 This, and a similar consideration for $\cdot_{\alpha}$, leads us
to define:
\begin{defn}

\begin{enumerate}
\item For $f,g\in\qc$ we define their quasi-sum $f\oplus g$ by 
\[
\left(f\oplus g\right)(x)=\sup_{y+z=x}\min\left\{ f(y),g(z)\right\} .
\]

\item For $f\in\qc$ and $\lambda>0$ we define $\lambda\odot f$ by 
\[
\left(\lambda\odot f\right)(x)=f\left(\frac{x}{\lambda}\right).
\]

\end{enumerate}
\end{defn}
For $\lambda=0$, we explicitly define 
\[
\left(0\odot f\right)(x)=\o_{\left\{ 0\right\} }(x)=\begin{cases}
1 & x=0\\
0 & x\ne0.
\end{cases}
\]
This definition ensures that $f\oplus(0\odot g)=f$ for every $f,g\in\qc$.
We use the notations $\oplus$ and $\odot$ instead of $\star_{-\infty}$
and $\cdot_{-\infty}$ because these operations will play a fundamental
role in the rest of this paper. 

So far we discussed properties of $\alpha$-concave functions which
made sense for every value of $\alpha$. We now want to state a few
results that are only true for quasi-concave functions and quasi-sums.
We will need the following definition:
\begin{defn}
For a function $f:\R^{n}\to[0,1]$ and $0<t\le1$ we define 
\[
K_{t}(f)=\left\{ x:\ f(x)\ge t\right\} 
\]
 to be the upper level sets of $f$. 
\end{defn}
We now have the following \foreignlanguage{american}{result, which
will play an important role in this note:}
\begin{thm}

\begin{enumerate}
\item Fix $f:\R^{n}\to[0,1]$. Then $f\in\qc$ if and only if $K_{t}(f)$
are compact, convex sets for all $0<t\le1$. 
\item For every $f,g\in\qc$,$\lambda>0$ and $0<t\le1$ we have 
\[
K_{t}\left(\left(\lambda\odot f\right)\oplus g\right)=\lambda K_{t}(f)+K_{t}(g).
\]

\end{enumerate}
\end{thm}
The sum $\oplus$ has another important property, we would now like
to discuss. Remember that if $f,g\in\Ca$, then $f,g\in\C{\alpha^{\prime}}$
for all $\alpha'<\alpha$, so we may look at the function $h=f\star_{\alpha^{\prime}}g\in\C{\alpha^{\prime}}$.
Generally, the function $h$ does not have to be in $\Ca$, even though
$f$ and $g$ are. Let us consider an example: choose $f(x)=g(x)=e^{-\left|x\right|}\in\lc$,
and choose $\alpha^{\prime}=-1$. In this case we have 
\[
\base{\left(-1\right)}f=\base{\left(-1\right)}g=e^{\left|x\right|}-1,
\]
 and since $h=f\star_{\alpha^{\prime}}g=2\cdot_{\alpha'}f$ we have
\[
\base{(-1)}h=2\cdot\left(e^{\frac{\left|x\right|}{2}}-1\right).
\]
Therefore 
\[
h(x)=\frac{1}{2e^{\frac{\left|x\right|}{2}}-1},
\]
and it is easy to check that $h\notin\lc$. Of course, we must have
$h\in\C{-1}$.

It turns out that such a situation cannot happen when $\alpha'=-\infty$:
\begin{thm}
If $f,g\in\C{\alpha}$ for some $-\infty\le\alpha\le0$, so does $f\oplus g$. 
\end{thm}
The proofs of the last two results will appear in \cite{Milman2012}.
Notice that by this theorem we have two different addition operations
on $\Ca$. One is $\star_{\alpha}$, and the second is the ``universal''
$\oplus$.

\section{Mixed integrals}

Recall the following theorem by Minkowski (see, e.g. \cite{Schneider1993}
for a proof):
\begin{thm*}
[Minkowski]Fix $K_{1},K_{2},\ldots,K_{m}\in\kc$. Then the function
$F:\left(\R^{+}\right)^{m}\to[0,\infty)$, defined by 
\[
F(\epsilon_{1},\epsilon_{2},\ldots,\epsilon_{m})=\vol\left(\epsilon_{1}K_{1}+\epsilon_{2}K_{2}+\cdots+\epsilon_{m}K_{m}\right),
\]
 is a homogenous polynomial of degree $n$, with non-negative coefficients.
\end{thm*}
The coefficients of this polynomial are called mixed volumes. To be
more exact, we have a function 
\[
V:\left(\kc\right)^{n}\to[0,\infty)
\]
 which is multilinear (with respect to the Minkowski sum), symmetric
(i.e. invariant to a permutation of its arguments), and which satisfies
$V(K,K,\ldots,K)=\vol(K)$. From these properties it is easy to deduce
that 
\[
F(\epsilon_{1},\epsilon_{2},\ldots,\epsilon_{m})=\sum_{i_{1},i_{2},\ldots,i_{n}=1}^{m}\epsilon_{i_{1}}\epsilon_{i_{2}}\cdots\epsilon_{i_{n}}\cdot V(K_{i_{1}},K_{i_{2}},\ldots,K_{i_{n}}).
\]
The number $V(K_{1},K_{2},\ldots,K_{n})$ is called the mixed volume
of the $K_{1},K_{2},\ldots,K_{n}$.

As we stated before, our goal is to prove a functional extension of
Minkowski's theorem, and to define a functional extension of mixed
volumes. We will state our results on $\qc$, since this is the largest
class of functions we consider, so all the results will be true for
every class $\Ca$ (and, in particular, on $\lc$). Of course, in
order to formulate and prove such a theorem, we need to decide what
are the functional analogs of volume, and of Minkowski sum.

For volume, if we want our theorem to be a true extension of Minkowski's,
we need a functional $\Phi$ on $\qc$ such that $\Phi\left(\o_{K}\right)=\vol\left(K\right)$.
A natural candidate is the Lebesgue integral, 
\[
\Phi(f)=\int_{\R^{n}}f(x)dx.
\]
 For the extension of addition, it turns out that the best possibility
is the quasi-sum $\oplus$. In fact, we have the following theorem:
\begin{thm}
\label{thm:func-Minkowski}Fix $f_{1},f_{2},\ldots,f_{m}\in\qc$.
Then the function $F:\left(\R^{+}\right)^{m}\to[0,\infty]$, defined
by 
\[
F(\epsilon_{1},\epsilon_{2},\ldots,\epsilon_{m})=\int\left[\left(\epsilon_{1}\odot f_{1}\right)\oplus\left(\epsilon_{2}\odot f_{2}\right)\oplus\cdots\oplus\left(\epsilon_{m}\odot f_{m}\right)\right]
\]
 is a homogenous polynomial of degree $n$, with non-negative coefficients.
\end{thm}
The proof will appear in \cite{Milman2012}. In complete analogy with
the case of convex bodies, this theorem is equivalent to the existence
of a function 
\[
V:\lc^{n}\to[0,\infty]
\]
 which is symmetric, multilinear (with respect to $\oplus$, of course)
and satisfies $V(f,f,\ldots,f)=\int_{\R^{n}}f(x)dx$. We will call
the number $V(f_{1},f_{2},\ldots f_{n})$ the mixed integral of $f_{1},f_{2},\ldots,f_{n}$.
In other words, the mixed integral is the polarization of the integral
$\int f$. We have the following representation formula for the mixed
integrals:
\begin{prop}
\label{prop:repr-formula}Fix $f_{1},f_{2},\ldots,f_{n}\in\qc$. Then
\[
V(f_{1},f_{2},\ldots,f_{n})=\int_{0}^{1}V\left(K_{t}(f_{1}),K_{t}(f_{2}),\ldots,K_{t}(f_{n})\right)dt.
\]

\end{prop}
Mixed integrals share many important properties with the classical
mixed volumes. We will mention a few in the following theorem:
\selectlanguage{american}%
\begin{thm}

\begin{enumerate}
\item For $K_{1},K_{2},\ldots,K_{n}\in\kc$ we have 
\[
V(K_{1},K_{2},\ldots,K_{n})=V(\o_{K_{1}},\o_{K_{2}},\ldots,\o_{K_{n}}).
\]

\item For every $f_{1},f_{2},\ldots,f_{n}\in\qc$ we have $V(f_{1},f_{2},\ldots,f_{n})\ge0$.
More generally, if we also have $g_{1},g_{2},\ldots,g_{n}\in\qc$
such that $f_{i}\ge g_{i}$ for all $i$, then 
\[
V(f_{1},f_{2},\ldots,f_{n})\ge V(g_{1},g_{2},\ldots,g_{n}).
\]
 
\item $V$ is rotation and translation invariant. Also, if we define 
\[
\left(uf\right)(x)=f(u^{-1}x)
\]
 for $f\in\qc$ and $u\in\GL(n)$, then 
\[
V(uf_{1},uf_{2},\ldots,uf_{n})=\left|\det u\right|\cdot V(f_{1},f_{2},\ldots,f_{n}).
\]
 
\selectlanguage{english}%
\item \label{enu:when-zero}Let $f_{1},f_{2},\ldots,f_{n}\in\qc$ and denote
by $K_{i}$ the support of $f_{i}$. Then $V(f_{1},f_{2},\ldots,f_{n})=0$
if and only if $V(K_{1},K_{2},\ldots,K_{n})=0$.
\item Fix an integer $1\le m\le n$ and functions $g_{m+1},\ldots,g_{n}\in\qc$.
Then the functional 
\[
\Phi(f)=V(f[m],g_{m+1},\ldots,g_{n})
\]
 satisfies a valuation type property: if $f_{1},f_{2}\in\qc$ and
$f_{1}\vee f_{2}=\max(f_{1},f_{2})\in\qc$ as well, then 
\[
\Phi\left(f_{1}\vee f_{2}\right)+\Phi\left(f_{1}\wedge f_{2}\right)=\Phi(f_{1})+\Phi(f_{2}).
\]
 Here $f_{1}\wedge f_{2}$ is an alternative notation for $\min\left\{ f_{1},f_{2}\right\} $. 
\end{enumerate}
\end{thm}
\selectlanguage{english}%
These properties are deduced by using Proposition \ref{prop:repr-formula}
and the corresponding properties for mixed volumes. We will prove
claim \ref{enu:when-zero}, and leave the others to the reader:
\begin{proof}
Denote $V=V(K_{1},K_{2},\ldots,K_{n})$, and define $f:(0,1]\to\R$
by 
\[
f(t)=V\left(K_{t}(f_{1}),K_{t}(f_{2}),\ldots,K_{t}(f_{n})\right).
\]
notice that $f$ is non-negative and non-increasing, by monotonicity
of mixed volumes. Since 
\[
K_{i}=\overline{\bigcup_{t>0}K_{t}(f_{i})},
\]
(the bar denotes the topological closure), and by continuity of mixed
volumes, we have $\lim_{t\to0^{+}}f(t)=V$.

Therefore, if $V=0$ then $f\equiv0$, so $V(f_{1},f_{2}\ldots,f_{n})=0$.
If, on the other hand, $V>0$, then $f(t)>\frac{V}{2}$ for all $t$
smaller than some $t_{0}$, so
\[
V(f_{1},f_{2},\ldots,f_{n})=\int_{0}^{1}f(t)dt\ge\int_{0}^{t_{0}}f(t)dt=\frac{t_{0}V}{2}>0.
\]

\end{proof}
A particularly interesting example of mixed volumes is quermassintegrals.
For $K\in\kc$, we define the $k$-th quermassintegral to be 
\[
W_{k}(K)=V(\underbrace{K,K,\ldots K}_{n-k\text{ times}},\underbrace{D,D,\ldots,D}_{k\text{ times}}),
\]
where $D$ is the Euclidean unit ball. Similarly, for $f\in\qc$ we
will define 
\[
W_{k}(f)=V(\underbrace{f,f,\ldots,f}_{n-k\text{ times}},\underbrace{\o_{D},\o_{D},\ldots,\o_{D}}_{k\text{ times}}).
\]
 By checking the definitions, we see that if $f\in\Ca$ and $K$ is
a convex set, then 
\[
\left[f\star_{\alpha'}\left(\epsilon\cdot_{\alpha'}\o_{K}\right)\right](x)=\sup_{y\in\epsilon K}f(x-y)
\]
 for every $\alpha'\le\alpha$. In particular, the left hand side
is independent of the exact value $\alpha'$. Since for $\alpha'=-\infty$
we obtain a polynomial in $\epsilon$, we must obtain the same polynomial
for every value of $\alpha'$. Therefore, as a direct corollary of
Theorem \ref{thm:func-Minkowski} we obtain the following statement,
which was independently obtained by Bobkov, Colesanti and Fragalà
(see \cite{Bobkov2012}):
\begin{prop}
\label{prop:func-Steiner}Fix $-\infty\le\alpha'\le\alpha\le0$. For
$f\in\Ca$ and $\epsilon>0$, define 
\[
f_{\epsilon}(x)=\left[f\star_{\alpha'}\left(\epsilon\cdot_{\alpha'}\o_{D}\right)\right](x)=\sup_{\left|y\right|\le\epsilon}f(x+y).
\]
 Then we have 
\[
\int f_{\epsilon}=\sum_{i=0}^{n}\binom{n}{i}W_{i}(f)\epsilon^{i}.
\]
 
\end{prop}
As stated, this result was also discovered by Bobkov, Colesanti and
Fragalà. Their paper continues to prove several properties of the
quermassintegrals, such as Prékopa-Leindler inequalities and a Cauchy-Kubota
integral formula. We will not pursue these points in this note. Let
us stress that Proposition \ref{prop:func-Steiner} only works for
quermassintegrals, where the different notions of sum happen to coincide.
For general mixed integrals, it is impossible to get polynomiality
for the operation $\star_{\alpha}$ unless $\alpha=-\infty$.

\section{Inequalities }

Now that we have a functional version of Minkowski's theorem, we would
like to prove inequalities between different mixed integrals. Let
us use the isoperimetric inequality as a test case. The classical
isoperimetric inequality, arguably the most famous inequality in geometry,
claims that for every (say convex) body $K\in\kc$ we have 
\[
S(K)\ge n\cdot\vol(D)^{\frac{1}{n}}\cdot\vol\left(K\right)^{\frac{n-1}{n}}.
\]
 Here $S(K)$ is the surface area of $K$, defined by 
\[
S(K)=\lim_{\epsilon\to0^{+}}\frac{\vol\left(K+\epsilon D\right)-\vol\left(K\right)}{\epsilon}=n\cdot W_{1}(K).
\]
 We would like to generalize this result to the functional setting.
The naive approach would be the try and bound $S(f):=nW_{1}(f)$ from
below using $\int f$. Unfortunately, this is impossible to do for
general quasi-concave functions. In fact, it is possible to construct
a sequence of functions $f_{k}\in\qc$ such that $\int f_{k}=1$ but
$S(f_{k})\to0$ as $k\to\infty$. We will present a concrete example
in \cite{Milman2012}.

Thus we will use a different approach, and prove an extension of the
isoperimetric inequality by recasting it as a rearrangement inequality.
In order to explain this idea, consider the following definition 
\begin{defn}
\label{def:SDR}
\begin{enumerate}
\item For a compact $K\in\kc$, define 
\[
K^{\ast}=\left(\frac{\vol(K)}{\vol(D)}\right)^{\frac{1}{n}}D.
\]
In other words, $K^{\ast}$ is the Euclidean ball with the same volume
as $K$. 
\item For $f\in\qc$, define its symmetric decreasing rearrangement $f^{\ast}$
using the relation 
\[
\KU_{t}\left(f^{\ast}\right)=\KU_{t}(f)^{\ast}.
\]

\end{enumerate}
\end{defn}
It is easy to see that this definition really defines a unique function
$f^{\ast}\in\qc$, which is rotation invariant.

Now, the isoperimetric inequality may be restated as $S(K)\ge S(K^{\ast})$
for $K\in\kc$. In this formulation, the functional extension turns
out to be true:
\begin{prop}
\label{prop:isoperimetric}If $f\in\qc$, then $S(f)\ge S(f^{\ast})$,
with equality if and only if $f$ is rotation invariant.
\end{prop}
This inequality is indeed an extension of the isoperimetric inequality,
as can be seen by choosing $f=\o_{K}$. It can also be useful for
general quasi-concave functions, because it reduces an $n$-dimensional
problem to a 1-dimensional one -- the function $f^{\ast}$ is rotation
invariant, and hence essentially ``one dimensional''. However, we
stress again that in general, this inequality does not yield a lower
bound for $S(f)$ in terms of $\int f$, as such a bound is impossible. 

Many other inequalities can be extended using similar formulations.
For example, the Brunn-Minkowski inequality states that for every
(say convex) sets $A,B\in\kc$ we have 
\[
\vol(A+B)^{\frac{1}{n}}\ge\vol(A)^{\frac{1}{n}}+\vol(B)^{\frac{1}{n}}.
\]
 Again, in general, it is impossible to bound $\int\left(f\oplus g\right)$
from below using $\int f$ and $\int g$. However, the Brunn-Minkowski
inequality may be written as $\left(A+B\right)^{\ast}\supseteq A^{\ast}+B^{\ast}$,
and in this representation it generalizes well:
\begin{thm}
\label{prop:BM}For every $f,g\in\qc$ we have $\left(f\oplus g\right)^{\ast}\ge f^{\ast}\oplus g^{\ast}.$ 
\end{thm}
In \cite{Milman2012} we will prove the above two results, as well
as extensions of the generalized Brunn-Minkowski inequalities for
mixed volumes and the Alexandrov-Fenchel inequality. We will not describe
these results here, since they require the notion of a ``generalized
rearrangement''. Instead, let us mention one elegant corollary of
the Alexandrov-Fenchel inequality. Remember that for convex bodies
we have the inequality 
\[
V(K_{1},K_{2},\ldots,K_{n})\ge\left[\prod_{i=1}^{n}\vol(K_{i})\right]^{\frac{1}{n}}.
\]
 The functional analog of this result is the following inequality:
\begin{thm}
\label{cor:from-AF}For all functions $f_{1},\ldots,f_{n}\in\qc$
we have 
\[
V(f_{1},f_{2},\ldots,f_{n})\ge V(f_{1}^{\ast},f_{2}^{\ast},\ldots,f_{n}^{\ast}).
\]
 
\end{thm}
Notice that Theorem \ref{cor:from-AF} is a generalization of the
isoperimetric inequality of Proposition \ref{prop:isoperimetric},
and it can also be used to deduce an Urysohn type inequality, bounding
$W_{n-1}(f)$ using $W_{n-1}(f^{\ast})$. 

If one is willing to restrict oneself to some class of $\alpha$-concave
functions, then it is suddenly possible to prove inequalities between
mixed integrals in a more familiar form. In order to state the result,
let us define for every $-\infty<\alpha\le0$ a function $g_{\alpha}\in\Ca$
by 
\[
g_{\alpha}(x)=\left(1+\frac{\left|x\right|}{\beta}\right)^{-\beta},
\]
 where, as usual $\beta=-\frac{1}{\alpha}$. Put differently, we choose
$g_{\alpha}$ to satisfy $\basea\left(g_{\alpha}\right)=\left|x\right|$
(for $\alpha=0$, we obtain $g_{0}(x)=e^{-\left|x\right|}$) . By
abuse of notation, we will also think of $g_{\alpha}$ as the function
from $[0,\infty)$ to $[0,\infty)$ defined by 
\[
g_{\alpha}(r)=\left(1+\frac{r}{\beta}\right)^{-\beta},
\]
 so $g_{\alpha}(x)=g_{\alpha}(\left|x\right|)$. We are now ready
to state:
\begin{thm}
\label{thm:alexandrov}Fix a function $f\in\Ca$ and integers $k$
and $m$ such that $0\le k<m<n$ . Then we have 
\[
\left(\frac{W_{k}(f)}{W_{k}(g_{\alpha})}\right)^{\frac{1}{n-k}}\le\left(\frac{W_{m}(f)}{W_{m}(g_{\alpha})}\right)^{\frac{1}{n-m}},
\]
 assuming $W_{k}(g_{\alpha})<\infty$. Equality occurs if and only
if $f=\lambda\odot g_{\alpha}$ for some $\lambda\ge0$. 
\end{thm}
Of course, if $W_{k}(g_{\alpha})=\infty$ the theorem is either trivial
(if $W_{k}(f)<\infty$) or meaningless (if $W_{k}(f)=\infty$). The
condition $W_{k}(g_{\alpha})<\infty$ is equivalent to $k>n+\frac{1}{\alpha}$,
and implies that all the other quantities in the theorem are finite
as well. Notice that since $k<m<n$ are all integers, we have $k\le n-2$.
Hence we need to choose $\alpha>-\frac{1}{2}$ for the theorem to
have any content.

In \cite{Milman2012}, a proof will be given for the case $\alpha=0$,
where the condition $W_{k}(g_{\alpha})<\infty$ is true for all $k$.
A key ingredient in the proof is a bound on the growth of moments
of log-concave functions. The general proof is similar, and depends
on the following lemma:
\begin{lem}
\label{lem:moments-comp}Let $f:[0,\infty)\to[0,\infty)$ be an $\alpha$-concave
function such that $f(0)=1$. Then for every $0\le k<m<-\frac{1}{\alpha}-1$
we have 
\[
\left(\frac{\int_{0}^{\infty}r^{m}f(r)dr}{\int_{0}^{\infty}r^{m}g_{\alpha}(r)dr}\right)^{\frac{1}{m+1}}\le\left(\frac{\int_{0}^{\infty}r^{k}f(r)dr}{\int_{0}^{\infty}r^{k}g_{\alpha}(r)dr}\right)^{\frac{1}{k+1}},
\]
 with equality if and only if 
\[
f(r)=\left(\lambda\odot g_{\alpha}\right)(r)=\left(1+\frac{r}{\lambda\beta}\right)^{-\beta}
\]
 for some $\lambda$. 
\end{lem}
Again, the condition $m<-\frac{1}{\alpha}-1$ simply ensures that
all of the integrals are finite. Under this condition the lemma follows
from Lemma 4.2 of \cite{Bobkov2009}, by taking $Q=g_{\alpha}$ (In
\cite{Bobkov2009} the equality condition is not explicitly stated,
but it can be deduced by carefully inspecting the proof).

We will conclude by sketching the proof of Theorem \ref{thm:alexandrov}.
Some parts of the proof, which are identical to the log-concave case,
will be glossed over and explained fully in \cite{Milman2012}.
\begin{proof}
First, we reduce the general case to the rotation invariant case,
by replacing $f$ with some generalized rearrangement $f^{W_{k}}$.
The definition of $f^{W_{k}}$ and its necessary properties will appear
in \cite{Milman2012}. 

So, assume without loss of generality that $f$ is rotation invariant.
By abuse of notation we will write $f(x)=f(\left|x\right|)$. A direct
computation (to also appear in \cite{Milman2012}) gives 
\[
W_{i}(f)=\left(n-i\right)\cdot\vol(D)\cdot\int_{0}^{\infty}r^{n-i-1}f(r)dr,
\]
 so 
\[
\left(\frac{W_{k}(f)}{W_{k}(g_{\alpha})}\right)^{\frac{1}{n-k}}=\left(\frac{\int_{0}^{\infty}r^{n-k-1}f(r)dr}{\int_{0}^{\infty}r^{n-k-1}g_{\alpha}(r)dr}\right)^{\frac{1}{n-k}},
\]
 and similarly for $m$. Remember that we assumed $W_{k}(g_{\alpha})<\infty$,
which is the same as
\[
\int_{0}^{\infty}r^{n-k-1}\left(1+\frac{r}{\beta}\right)^{-\beta}dr<\infty.
\]
 This implies that $n-k-1-\beta<-1$, or $k>n+\frac{1}{\alpha}$,
like we claimed. In particular we have 
\[
0\le n-m-1<n-k-1<-\frac{1}{\alpha}-1,
\]
so we can use Lemma \ref{lem:moments-comp} and conclude that
\[
\left(\frac{\int_{0}^{\infty}r^{n-k-1}f(r)dr}{\int_{0}^{\infty}r^{n-k-1}g_{\alpha}(r)dr}\right)^{\frac{1}{n-k}}\le\left(\frac{\int_{0}^{\infty}r^{n-m-1}f(r)dr}{\int_{0}^{\infty}r^{n-m-1}g_{\alpha}(r)dr}\right)^{\frac{1}{n-m}}.
\]
 This is the same as 
\[
\left(\frac{W_{k}(f)}{W_{k}(g_{\alpha})}\right)^{\frac{1}{n-k}}\le\left(\frac{W_{m}(f)}{W_{m}(g_{\alpha})}\right)^{\frac{1}{n-m}},
\]
which is what we wanted. 

The equality case will follow from the equality case of Lemma \ref{lem:moments-comp},
but we will not give the details here.
\end{proof}
\bibliographystyle{plain}
\bibliography{hyperbolic}

\end{document}